\documentclass[a4paper,11pt]{amsart}
\usepackage{graphicx}
\usepackage{mathptmx}
\usepackage{amsmath}
\usepackage{amssymb}
\usepackage{enumitem}
\usepackage{xcolor}
\usepackage{xparse}
\NewDocumentCommand{\eulerian}{omm}
 {%
  \genfrac<>{0pt}{}{#2}{#3}%
  \IfValueT{#1}{_{\!#1}}%
 }

 \newmuskip\pFqmuskip

\newcommand*\pFq[6][8]{%
  \begingroup 
  \pFqmuskip=#1mu\relax
  \mathchardef\normalcomma=\mathcode`,
  \mathcode`\,=\string"8000
  \begingroup\lccode`\~=`\,
  \lowercase{\endgroup\let~}\pFqcomma
  {}_{#2}F_{#3}{\left(\genfrac..{0pt}{}{#4}{#5}\bigg|#6\right)}%
  \endgroup
}
\newcommand{\pFqcomma}{{\normalcomma}\mskip\pFqmuskip}

\newtheorem{theorem}{Theorem}
\newtheorem{lemma}[theorem]{Lemma}

\begin{document}

\title[Multi-Lah numbers and multi-Stirling numbers of the first kind]{Multi-Lah numbers and multi-Stirling numbers of the first kind}

\author{Dae San KiM}
\address{Department of Mathematics, Sogang University, Seoul 121-742, Republic of Korea}
\email{dskim@sogang.ac.kr}

\author{Hye Kyung KiM}
\address{Department of Mathematical Education, Daegu Catholic University,
	Gyeongsan 38430, Republic of Korea}
\email{hkkim@cu.ac.kr}

\author{Taekyun  Kim}
\address{Department of Mathematics, Kwangwoon University, Seoul 139-701, Republic of Korea}
\email{taekyun64@hotmail.com}

\author{Hyunseok  Lee}
\address{Department of Mathematics, Kwangwoon University, Seoul 139-701, Republic of Korea}
\email{luciasconstant@kw.ac.kr}

\author{SeongHo Park}
\address{Department of Mathematics, Kwangwoon University, Seoul 139-701, Republic of Korea}
\email{abcd2938471@kw.ac.kr}


\begin{abstract}
In this paper, we introduce multi-Lah numbers and multi-Stirling numbers of the first kind and recall multi-Bernoulli numbers, all of whose generating functions are given with the help of multiple logarithm. 
The aim of this paper is to study several relations among those three numbers. In more detail, we represent the multi-Bernoulli numbers in terms of the multi-Stirling numbers of the first kind and vice versa, and the multi-Lah numbers in terms of multi-Stirling numbers. In addition, we deduce a recurrence relation for multi-Lah numbers.
\end{abstract}

\maketitle

\section{Introduction}

As is well known, the unsigned Stirling number ${n\brack r}$ counts the number of permutations of a set with $n$ elements which are products of $r$ disjoint cycles. We generalize these numbers to the multi-Stirling numbers of the first kind $S_{1}^{(k_{1},k_{2},\dots,k_{r})}(n,r)$ (see \eqref{10}) which reduce to the unsigned Stirling numbers of the first kind for $(k_{1},k_{2},\dots,k_{r})=(1,1,\dots,1)$. Indeed, $S_{1}^{(1,1,\dots,1)}(n,r)={n\brack r}$. \par
It is also well known that the unsigned Lah numbers $L(n,k)$ counts the number of ways of a set of $n$ elements can be partitioned into $k$ nonempty linearly ordered subsets. These numbers are generalized to the multi-Lah numbers $L^{(k_{1},k_{2},\dots,k_{r})}(n,r)$ (see \eqref{18}) which reduce to the unsigned Lah numbers for $(k_{1},k_{2},\dots,k_{r})=(1,1,\dots,1)$. In fact, $L^{(1,1,\dots,1)}(n,r)=L(n,r).$ \par
In addition, we need to recall the multi-Bernoulli numbers $B_{n}^{(k_1,k_2,\dots,k_{r})}$ (see \eqref{9}) which were introduced earlier under the different name the generalized Bernoulli numbers with order $r$ in [7]. These numbers reduce to the Bernoulli numbers of order $r$ up to some constants. Indeed, we see that $B_{m}^{(1,1,\dots,1)}=\frac{1}{r!}(-1)^{m}B_{m}^{(r)}$. \par
The common feature of those three kinds of numbers is that they are all defined with the help of the multiple logarithm $\mathrm{Li}_{k_{1},k_{2},\dots,k_{r}}(z)$ (see \eqref{8}), which reduce to the poly-logarithm $\mathrm{Li}_{k_{1}}(z)$, for $r=1$. \par
The aim of this paper is to study several relations among those three kinds of numbers. In more detail, we represent the multi-Bernoulli numbers in terms of the multi-Stirling numbers of the first kind and vice versa, and the multi-Lah numbers in terms of multi-Stirling numbers. Moreover, we deduce a recurrence relation for multi-Lah numbers.  For the rest of this section, we recall the necessary facts that will be needed throughout this paper.

\vspace{0.1in}

The unsigned Stirling numbers $L(n,k)$ are defined by 
\begin{equation}
\langle x\rangle_{n}=\sum_{k=0}^{n}L(n,k)(x)_{k},\quad (n\ge 0),\quad (\mathrm{see}\ [2,3,4,5,6]),\label{1}
\end{equation}
where $\langle x\rangle_{0}=1,\ \langle x\rangle_{n}=x(x+1)\cdots(x+n-1),\ (n \ge 1), $  and  $(x)_{0}=1,\ (x)_{n}=x(x-1)\cdots(x-n+1), \ (n \ge 1)$. \par 
The inverse formula of \eqref{1} is given by 
\begin{displaymath}
(x)_{n}=\sum_{k=0}^{n}(-1)^{n-k}L(n,k)\langle x\rangle_{k},
\quad(n\ge 0).	
\end{displaymath}
From \eqref{1}, we can derive the generating function of unsigned Lah numbers given by 
\begin{equation}
\frac{1}{k!}\bigg(\frac{t}{1-t}\bigg)^{k}	=\sum_{n=k}^{\infty}L(n,k)\frac{t^{n}}{n!},\quad (\mathrm{see}\ [2,3,8]).\label{2}
\end{equation}
Thus, we note that 
\begin{equation}
L(n,k)=\binom{n-1}{k-1}\frac{n!}{k!},\quad (n,k\ge 1),\quad(\mathrm{see}\ [2,8]).\label{3}
\end{equation}
The Stirling numbers of the first kind are defined by 
\begin{equation}
(x)_{n}=\sum_{k=0}^{n}S_{1}(n,k)x^{k},\quad (n\ge 0),\quad (\mathrm{see}\ [5,8]),\label{4}	
\end{equation}
and the Stirling numbers of the second kind are defined by 
\begin{equation}
x^{n}=\sum_{k=0}^{n}S_{2}(n,k)(x)_{k},\quad (n\ge 0),\quad (\mathrm{see}\ [2,3,8]).\label{5}
\end{equation}
From \eqref{4} and \eqref{5}, we note that 
\begin{equation}
\frac{1}{k!}\big(e^{t}-1\big)^{k}=\sum_{n=k}^{\infty}S_{2}(n,k)\frac{t^{n}}{n!},\quad (\mathrm{see}\ [2,4,6,8]),\label{6}	
\end{equation}
and 
\begin{equation}
\frac{1}{k!}\big(\log(1+t)\big)^{k}=\sum_{n=k}^{\infty}S_{1}(n,k)\frac{t^{n}}{n!},\quad (\mathrm{see}\ [8]).\label{7}	
\end{equation}
For any $k_{i}\ge 1\ (1\le i\le r)$,\ and $|z|<1$, the multiple logarithm is defined by 
\begin{equation}
\mathrm{Li}_{k_{1},k_{2},\dots,k_{r}}(z)=\sum_{0<m_{1}<m_{2}<\cdots<m_{r}}\frac{z^{m_{r}}}{m_{1}^{k_{1}}m_{2}^{k_{2}}\cdots m_{r}^{k_{r}}},\quad(\mathrm{see}\ [1,7]).\label{8} 
\end{equation}
If $r=1$, \,$\displaystyle \mathrm{Li}_{k_{1}}(z)=\sum_{m=1}^{\infty}\frac{z^m}{m^{k_{1}}}\displaystyle$ is the poly-logarithm. \par 
The multi-Bernoulli numbers, which are called the generalized Bernoulli numbers with order $r$ in [7], are defined by 
\begin{equation}
\frac{\mathrm{Li}_{k_{1},k_{2},\dots,k_{r}}(z)}{z^{r}}\bigg|_{z=1-e^{-t}}=\sum_{m=0}^{\infty}\frac{B_{m}^{(k_{1},k_{2},\dots,k_{r})}}{m!}t^{m}. \label{9}
\end{equation}
From \eqref{13}, we note that 
\begin{displaymath}
B_{m}^{\overbrace{(1,1,\dots,1)}^{r-times}}=\frac{1}{r!}(-1)^{m}B_{m}^{(r)},\quad (m\ge 0,\ r\ge 1),\quad(\mathrm{see}\ [7]),
\end{displaymath}
where $B_{n}^{(r)}$ are the Bernoulli numbers of order $r$ given by  
\begin{displaymath}
\bigg(\frac{t}{e^{t}-1}\bigg)^{r}=\sum_{n=0}^{\infty}B_{n}^{(r)}\frac{t^{n}}{n!},\quad (\mathrm{see}\ [1,7,8]).
\end{displaymath}

\section{Multi-Lah numbers and multi-Stirling numbers of the first kind} 
Now, we define the {\it{multi-Stirling numbers of the first kind}} by
\begin{equation}
\mathrm{Li}_{k_{1},k_{2},\dots,k_{r}}(t)=\sum_{n=r}^{\infty}S_{1}^{(k_{1},k_{2},\dots,k_{r})}(n,r)\frac{t^{n}}{n!}, \label{10}	
\end{equation}
where $k_{i}\ge 1\ (1\le i\le r-1)$, $k_{r}\ge 2, \ and\ |t|<1$. \par 
From \eqref{8}, we note that 
\begin{align}
\frac{d}{dt} \mathrm{Li}_{k_{1},k_{2},\dots,k_{r}}(t)\ &=\ \frac{d}{dt}\sum_{0<m_{1}<m_{2}<\cdots<m_{r}}\frac{t^{m_{r}}}{m_{1}^{k_{1}}m_{2}^{k_{2}}\cdots m_{r}^{k_{r}}}\label{11} \\
&=\ \frac{1}{t}\sum_{0<m_{1}<m_{2}<\cdots<m_{r}}\frac{t^{m_{r}}}{m_{1}^{k_{1}}m_{2}^{k_{2}}\cdots m_{r-1}^{k_{r-1}} m_{r}^{k_{r}-1}}\nonumber\\
&=\ \frac{1}{t} \mathrm{Li}_{k_{1},k_{2},\dots,k_{r-1},k_{r}-1}(t).\nonumber
\end{align}
Let us take $k_{r}=1$ in \eqref{11}. Then we have 
\begin{align}
\frac{d}{dt} \mathrm{Li}_{k_{1},k_{2},\dots,k_{r-1},1}(t)\ &=\ \frac{1}{t}\mathrm{Li}_{k_{1},k_{2},\dots,k_{r-1},0}(t)\label{12} \\
&=\ \frac{1}{t}\sum_{0<m_{1}<\cdots<m_{r-1}}\frac{1}{m_{1}^{k_{1}}\cdots m_{r-1}^{k_{r-1}}}\sum_{m_{r}=m_{r-1}+1}^{\infty}t^{m_{r}}	\nonumber \\
&=\ \sum_{0<m_{1}<m_{2}<\cdots<m_{r-1}}\frac{1}{m_{1}^{k_{1}}m_{2}^{k_{2}}\cdots m_{r-1}^{k_{r-1}}}\frac{t^{m_{r-1}+1}}{1-t}\frac{1}{t}\nonumber \\
&=\ \frac{1}{1-t} \mathrm{Li}_{k_{1},k_{2},\dots,k_{r-1}}(t).\nonumber
\end{align}
We claim that the following relations hold. For this, we only need to show the first equality which we prove by induction on $r$.
\begin{equation}
\mathrm{Li}_{\underbrace{1,1,\dots,1}_{r-times}}(t)=\frac{1}{r!}\big(-\log(1-t)\big)^{r}=\sum_{n=r}^{\infty}(-1)^{n-r}S_{1}(n,r)\frac{t^{n}}{n!}.\label{13}
\end{equation}
If $r=1$, then $\mathrm{Li}_{1}(t)=\sum_{m=1}^{\infty}\frac{t^m}{m}=-\log(1-t)$, as we wanted.
Assume that $r \ge 2$ and that it holds for $r-1$. By \eqref{12} and induction hypothesis, we get
\begin{equation}
\frac{d}{dt} \mathrm{Li}_{\underbrace{1,1,\dots,1}_{r-times}}(t) = \frac{1}{1-t} \mathrm{Li}_{\underbrace{1,1,\dots,1}_{(r-1)-times}}(t)=\frac{1}{(r-1)!}\frac{1}{1-t}\big(-\log (1-t)\big)^{r-1}.\label{14}
\end{equation}
Now, by \eqref{14} we obtain
\begin{align}
\mathrm{Li}_{\underbrace{1,1,\dots,1}_{r-times}}(t)& =\frac{1}{(r-1)!}\int_{0}^{t}\frac{1}{1-t}\big(-\log (1-t)\big)^{r-1} dt  \label{15} \\
&=\frac{1}{(r-1)!}\int_{0}^{-\log (1-t)}u^{r-1} du=\frac{1}{r!}\big(-\log(1-t)\big)^{r}. \nonumber
\end{align}
Thus our proof is completed.
From \eqref{10}, we note that 
\begin{equation}
\mathrm{Li}_{\underbrace{1,1,\dots,1}_{r-times}}(t)=\sum_{n=r}^{\infty}S_{1}^{\overbrace{(1,1,\dots,1)}^{r-times}}(n,r)\frac{t^{n}}{n!}.\label{16}	
\end{equation}
Therefore, by \eqref{13} and \eqref{16}, we obtain the following lemma. 
\begin{lemma} For 
$n,r\ge 1$, we have 	
\begin{displaymath}
S_{1}^{\overbrace{(1,1,\dots,1)}^{r-times}}(n,r)=(-1)^{n-r}S_{1}(n,r)={n\brack r},
\end{displaymath}
where ${n\brack r}$ are the unsigned Stirling numbers of the first kind. 
\end{lemma}
We observe that 
\begin{align}
\frac{\mathrm{Li}_{k_{1},k_{2},\dots,k_{r}}(1-e^{-t})}{(1-e^{-t})^{r}}\ &=\ \frac{1}{(1-e^{-t})^{r}}\sum_{m=r}^{\infty}S_{1}^{(k_{1},k_{2},\dots,k_{r})}(m,r)\frac{1}{m!}\big(1-e^{-t}\big)^{m} \label{17}\\
&=\ \frac{1}{(1-e^{-t})^{r}}\sum_{l=r}^{\infty}\sum_{m=r}^{l}(-1)^{m-l}S_{2}(l,m)S_{1}^{(k_{1},k_{2},\dots,k_{r})}(m,r)\frac{t^{l}}{l!}\nonumber \\
&=\ \frac{(-1)^{r}t^{r}}{(1-e^{-t})^{r}}\sum_{l=0}^{\infty}\sum_{m=r}^{l+r}(-1)^{m-l}S_{2}(l+r,m)S_{1}^{(k_{1},\dots,k_{r})}(m,r)\frac{l!}{(l+r)!}\frac{t^{l}}{l!}\nonumber\\
&=\ \sum_{j=0}^{\infty}B_{j}^{(r)}(-1)^{r-j}\frac{t^{j}}{j!}\sum_{l=0}^{\infty}\bigg(\sum_{m=r}^{l+r}(-1)^{m-l}S_{2}(l+r,m)S_{1}^{(k_{1},\dots,k_{r})}(m,r) \frac{1}{r!\binom{l+r}{l}}\bigg)\frac{t^{l}}{l!}\nonumber \\
&=\ \sum_{n=0}^{\infty}\bigg(\sum_{l=0}^{n}\sum_{m=r}^{l+r}\frac{\binom{n}{l}B_{n-l}^{(r)}(-1)^{n-r-m}}{r!\binom{l+r}{r}}S_{2}(l+r,m)S_{1}^{(k_{1},\dots,k_{r})}(m,r)\bigg)\frac{t^{n}}{n!}.\nonumber	
\end{align}
Therefore, by \eqref{9} and \eqref{17}, we obtain the following theorem. 
\begin{theorem}
For $k_{i}\ge 1\ (i=1,2,\dots,r),\ and \ n\ge 0$, we have 
\begin{displaymath}
B_{n}^{(k_{1},k_{2},\dots,k_{r})}= \sum_{l=0}^{n}\sum_{m=r}^{l+r}\frac{\binom{n}{l}B_{n-l}^{(r)}(-1)^{n-r-m}}{r!\binom{l+r}{r}}S_{2}(l+r,m)S_{1}^{(k_{1},\dots,k_{r})}(m,r)	.
\end{displaymath}
\end{theorem}
For any integer $k_{i}\ (i=1,2,\dots,r)$, in the view of \eqref{9}, we define $L^{(k_{1},k_{2},\dots,k_{r})}(n,r), (n,r\ge 0)$, which are called {\it{multi-Lah numbers}}, as
\begin{equation}
\frac{\mathrm{Li}_{k_{1},k_{2},\dots,k_{r}}(1-e^{-t})}{(1-t)^{r}}=\sum_{n=r}^{\infty} L^{(k_{1},k_{2},\dots,k_{r})}(n,r)\frac{t^{n}}{n!}.\label{18}
\end{equation}
From \eqref{13}, we note that 
\begin{align}
\sum_{n=r}^{\infty}L^{\overbrace{(1,1,\dots,1)}^{r-times}}(n,r)\frac{t^{n}}{n!}\ &=\ \frac{1}{(1-t)^{r}}\mathrm{Li}_{\underbrace{1,1,\dots,1}_{r-times}}\big(1-e^{-t}\big)\label{19} \\
&=\ \frac{1}{(1-t)^{r}}\frac{1}{r!}\Big(-\log\big(1-(1-e^{-t})\big)\Big)^{r}	\nonumber \\
&=\ \frac{1}{r!}\bigg(\frac{t}{1-t}\bigg)^{r}\ =\ \sum_{n=r}^{\infty}L(n,r)\frac{t^{n}}{n!}.\nonumber
\end{align}
Thus, by \eqref{19}, we get 
\begin{displaymath}
L^{\overbrace{(1,1,\dots,1)}^{r-times}}(n,r)=L(n,r),\quad (n,r\ge 0).
\end{displaymath}
For $n\ge 1$, from \eqref{8} and \eqref{18}, we get 
\begin{align}
&\sum_{n=r}^{\infty} L^{(k_{1},\dots,k_{r-1},-k_{r})}(n,r)\frac{t^{n}}{n!}=\frac{1}{(1-t)^{r}} \mathrm{Li}_{k_{1},\dots,k_{r-1},-k_{r}}(1-e^{t}) \label{20}\\
&=\frac{1}{(1-t)^{r}}\sum_{0<m_{1}<m_{2}<\cdots<m_{r-1}}\frac{1}{m_{1}^{k_{1}}m_{2}^{k_{2}}\cdots m_{r-1}^{k_{r-1}}}\sum_{m_{r}=m_{r-1}+1}^{\infty}\frac{(1-e^{-t})^{m_{r}}}{m_{r}^{-k_{r}}}\nonumber\\
&=\frac{1}{(1-t)^{r}}\sum_{0<m_{1}<m_{2}<\cdots<m_{r-1}}\frac{1}{m_{1}^{k_{1}}m_{2}^{k_{2}}\cdots m_{r-1}^{k_{r-1}}}\sum_{m_{r}=1}^{\infty}\frac{(1-e^{-t})^{m_{r}+m_{r-1}}}{(m_{r}+m_{r-1})^{-k_{r}}}\nonumber\\
&=\frac{1}{(1-t)^{r}}\sum_{0<m_{1}<m_{2}<\cdots<m_{r-1}}\frac{(1-e^{t})^{m_{r-1}}}{m_{1}^{k_{1}}m_{2}^{k_{2}}\cdots m_{r-1}^{k_{r-1}}}\sum_{m_{r}=1}^{\infty}\frac{(-1)^{m_{r}}m_{r}!}{(m_{r}+m_{r-1})^{-k_{r}}}\frac{1}{m_{r}!}\big(e^{-t}-1\big)^{m_{r}}\nonumber\\ 
&=\frac{1}{(1-t)^{r}}\sum_{0<m_{1}<m_{2}<\cdots<m_{r-1}}\frac{(1-e^{t})^{m_{r-1}}}{m_{1}^{k_{1}}m_{2}^{k_{2}}\cdots m_{r-1}^{k_{r-1}}}\sum_{m_{r}=1}^{\infty}\frac{(-1)^{m_{r}}m_{r}!}{(m_{r}+m_{r-1})^{-k_{r}}}\sum_{l=m_{r}}^{\infty}S_{2}(l,m_{r})\frac{(-t)^{l}}{l!}\nonumber\\
 &=\frac{1}{(1-t)^{r}}\sum_{0<m_{1}<m_{2}<\cdots<m_{r-1}}\frac{(1-e^{t})^{m_{r-1}}}{m_{1}^{k_{1}}m_{2}^{k_{2}}\cdots m_{r-1}^{k_{r-1}}}\sum_{l=1}^{\infty}\bigg(\sum_{m_{r}=1}^{l}\frac{(-1)^{m_{r}+l}m_{r}!}{(m_{r}+m_{r-1})^{-k_{r}}}S_{2}(l,m_{r})\bigg)\frac{t^{l}}{l!}\nonumber \\
 &=\frac{1}{1-t}\sum_{l=1}^{\infty}\bigg(\sum_{m_{r}=1}^{l}(-1)^{m_{r}+l}m_{r}!S_{2}(l,m_{r})\sum_{j=0}^{k_{r}}\binom{k_{r}}{j}m_{r}^{k_{r}-j} \nonumber \\
 &\quad\quad\times \bigg(\frac{1}{1-t}\bigg)^{r-1}\sum_{0<m_{1}<m_{2}<\cdots<m_{r-1}}\frac{(1-e^{-t})^{m_{r-1}}}{m_{1}^{k_{1}}\cdots m_{r-1}^{k_{r-1}-j}}\bigg)\frac{t^{l}}{l!}\nonumber \\
 &=\frac{1}{1-t}\sum_{l=1}^{\infty}\bigg(\sum_{m_{r}=1}^{l}(-1)^{m_{r}+l}m_{r}!S_{2}(l,m_{r})\sum_{j=0}^{k_{r}}\binom{k_{r}}{j}m_{r}^{k_{r}-j}\bigg(\frac{1}{1-t}\bigg)^{r-1}\mathrm{Li}_{k_{1},\dots,k_{r-1}-j}(1-e^{-t})\bigg)\frac{t^{l}}{l!}.\nonumber 
\end{align}
From \eqref{20}, we note that 
\begin{align}
&\sum_{n=r}^{\infty}L^{(k_{1},k_{2},\dots,k_{r-1},-k_{r})}(n,r)\frac{t^{n}}{n!} \label{21}\\
&\ =\frac{1}{1-t}\sum_{l=1}^{\infty}\bigg(\sum_{m_{r}=1}^{l}\sum_{j=0}^{k_{r}}(-1)^{m_{r}+l}m_{r}!S_{2}(l,m_{r})\binom{k_{r}}{j}m_{r}^{k_{r}-j}\nonumber \\
&\qquad \times \sum_{m=r-1}^{\infty}L^{(k_{1},k_{2},\dots,k_{r-2},k_{r-1}-j)}(m,r-1)\frac{t^{m}}{m!}\bigg)\frac{t^{l}}{l!}	\nonumber \\
&\ =\frac{1}{1-t}\sum_{p=r}^{\infty}\bigg(\sum_{l=1}^{p+1-r}\binom{p}{l}\sum_{m_{r}=1}^{p}\sum_{j=0}^{k_{r}}(-1)^{m_{r}+l}m_{r}!m_{r}^{k_{r}-j}S_{2}(l,m_{r})\binom{k_{r}}{j}\nonumber \\
&\qquad \times L^{(k_{1},k_{2},\dots,k_{r-2},k_{r-1}-j)}(p-l,r-1)\bigg)\frac{t^{p}}{p!} \nonumber \\
&= \sum_{n=r}^{\infty}\bigg(\sum_{p=r}^{n} \sum_{l=1}^{p+1-r}\binom{p}{l}\sum_{m_{r}=1}^{p}\sum_{j=0}^{k_{r}}(-1)^{m_{r}+l}m_{r}!m_{r}^{k_{r}-j}S_{2}(l,m_{r})\binom{k_{r}}{j}\nonumber \\
&\qquad\times L^{(k_{1},k_{2},\dots,k_{r-2},k_{r-1}-j)}(p-l,r-1)\frac{n!}{p!}\bigg)\frac{t^{n}}{n!}.\nonumber
\end{align}
Therefore, by comparing the coefficients on both sides of \eqref{21}, we obtain the following theorem. 
\begin{theorem}
For any $k_{i}\ge 1\ (i=1,2,\dots,r)$, and $n,r\in\mathbb{N}$, we have 
\begin{align*}
& L^{(k_{1},k_{2},\dots,k_{r-1},-k_{r})}(n,r)\\
&= \sum_{p=r}^{n} \sum_{l=1}^{p}\sum_{m_{r}=1}^{p}\sum_{j=0}^{k_{r}}(-1)^{m_{r}+l}m_{r}!m_{r}^{k_{r}-j}\binom{p}{l}\binom{k_{r}}{j}S_{2}(l,m_{r}) L^{(k_{1},k_{2},\dots,k_{r-2},k_{r-1}-j)}(p-l,r-1)\frac{n!}{p!}.
\end{align*}
\end{theorem}
Replacing $t$ by $-\log(1-t)$ in \eqref{9}, we get 
\begin{align}
\frac{1}{t^{r}}\mathrm{Li}_{k_{1},k_{2},\dots,k_{r}}(t)\ &=\ \sum_{m=0}^{\infty}B_{m}^{(k_{1},k_{2},\dots,k_{r})}(-1)^{m}\frac{1}{m!}\big(\log(1-t)\big)^{m}\nonumber \\
&=\ \sum_{m=0}^{\infty}B_{m}^{(k_{1},k_{2},\dots,k_{r})}(-1)^{m}\sum_{n=m}^{\infty}S_{1}(n,m)(-1)^{n}\frac{t^{n}}{n!}\label{22} \\
&=\ \sum_{n=0}^{\infty}\bigg(\sum_{m=0}^{n}(-1)^{n-m} B_{m}^{(k_{1},k_{2},\dots,k_{r})}S_{1}(n,m)\bigg)\frac{t^{n}}{n!}.\nonumber	
\end{align}
On the other hand, by \eqref{10}, we get 
\begin{align}
\frac{1}{t^{r}} \mathrm{Li}_{k_{1},k_{2},\dots,k_{r}}(t)\ &=\ \frac{1}{t^{r}}\sum_{n=r}^{\infty}S_{1}^{(k_{1},k_{2},\dots,k_{r})}(n,r)\frac{t^{n}}{n!} \label{23} \\
&=\ \sum_{n=0}^{\infty}S_{1}^{(k_{1},k_{2},\dots,k_{r})}(n+r,r)\frac{n!}{(n+r)!}\frac{t^{n}}{n!}\nonumber \\
&=\ \sum_{n=0}^{\infty}\frac{S_{1}^{(k_{1},k_{2},\dots,k_{r})}(n+r,r)}{r!\binom{n+r}{n}}\frac{t^{n}}{n!}. \nonumber 	
\end{align}
Therefore, by \eqref{22} and \eqref{23}, we obtain the following theorem. 
\begin{theorem}
For each $k_{i}\ge 1\ (i=1,2,\dots,r),\ n\ge 0,\,\,and\  r\in\mathbb{N}$, we have 
\begin{displaymath}
S_{1}^{(k_{1},k_{2},\dots,k_{r})}(n+r,r)=r!\binom{n+r}{n}\sum_{m=0}^{n}B_{m}^{(k_{1},k_{2},\dots,k_{r})}{n \brack m}.
\end{displaymath}	
\end{theorem}
Now, we observe that 
\begin{align}
&\sum_{n=r}^{\infty}L^{(k_{1},k_{2},\dots,k_{r})}(n,r)\frac{t^{n}}{n!}=\frac{1}{(1-t)^{r}}\mathrm{Li}_{k_{1},k_{2},\dots,k_{r}}(1-e^{-t})\label{24} \\
&\ =\sum_{j=0}^{\infty}\binom{r+j-1}{j}t^{j}\sum_{m=r}^{\infty}S_{1}^{k_{1},k_{2},\dots,k_{r}}(m,r)\frac{1}{m!}\big(1-e^{-t}\big)^{m} \nonumber \\ 
&\ =\sum_{j=0}^{\infty}\binom{r+j-1}{j}t^{j}\sum_{l=r}^{\infty}\bigg(\sum_{m=r}^{l}S_{1}^{(k_{1},\dots,k_{r})}(m,r)(-1)^{m-l}S_{2}(l,m)\bigg)\frac{t^{l}}{l!} \nonumber \\
&\ =\sum_{n=r}^{\infty}\bigg(\sum_{l=r}^{n}\sum_{m=r}^{l}S_{1}^{(k_{1},\dots,k_{r})}(m,r)(-1)^{m-l}S_{2}(l,m)\binom{r+n-l-1}{n-l}\frac{n!}{l!}\bigg)\frac{t^{n}}{n!}.\nonumber
\end{align}
Therefore, by comparing the coefficients on both sides of \eqref{24}, we obtain the following theorem. 
\begin{theorem}
For each $k_{i}\ (i=1,2,\dots,r),\,\,and\,\, n,r\in\mathbb{N}$, with $n\ge r$, we have 
\begin{align*}
	& L^{(k_{1},k_{2},\dots,k_{r})}(n,r)\\
	&\quad= \sum_{l=r}^{n}\sum_{m=r}^{l}S_{1}^{(k_{1},\dots,k_{r})}(m,r)(-1)^{m-l}S_{2}(l,m)\binom{r+n-l-1}{n-l}\frac{n!}{l!}.
\end{align*}	
\end{theorem}

\section{Conclusion} 
There are various ways of studying special polynomials and numbers which include generating functions, combinatorial methods, $p$-adic analysis, umbral calculus, special functions, differential equations and probability theory. In this paper, using generating function method and by making use of the multiple logarithm, we studied three numbers, namely the multi-Stirling numbers of the first kind, the multi-Lah numbers and the multi-Bernoulli numbers which reduce respectively to the unsigned Stirling numbers of the first kind, the Lah numbers and the higher-order Bernoulli numbers up to constants when the index is specialized to $(k_{1},k_{2},\dots,k_{r})=(1,1,\dots,1)$. We deduced several relations among those numbers. In more detail, we expressed the multi-Bernoulli numbers in terms of the multi-Stirling numbers of the first kind and vice versa, and the multi-Lah numbers in terms of multi-Stirling numbers. Further,  we derived a recurrence relation for multi-Lah numbers. \par
It is our continuous interest to explore some special numbers and polynomials by using different tools like the ones mentioned in the above.

\vspace{0.2in}

\noindent{\bf{Acknowledgments:}} 
Not applicable.
\vspace{0.2in}

\noindent{\bf{Funding:}} 
This work was supported by the Basic Science Research Program, the National  Research Foundation of Korea, (NRF-2021R1F1A1050151).

\vspace{0.2in}

\noindent{\bf {Availability of data and materials:}}
Not applicable.

\vspace{0.1in}

\noindent{\bf {Competing interests:}}
The authors declare no conflict of interest.

\vspace{0.1in}

\noindent{\bf{Authors' contributions:}}  D.S.K., H.K.K and T.K. conceived of the framework and structured the whole paper; T. K. and D.S.K. wrote the paper; D.S.K., H.K.K. and T.K. completed the revision of the article; H.L. and S.P checked the errors of the article.. All authors have read and agreed to the published version of the manuscript.

\vspace{0.1in}



\vspace{0.1in}



\begin{thebibliography}{9}
\bibitem{1}
Borwein, J. M.; Chan, O-Y. \emph{Duality in tails of multiple-zeta values,} Int. J. Number Theory \textbf{6} (2010), no. 3, 501-514. 
\bibitem{2}
Comtet, L. \emph{Advanced combinatorics,} The art of finite and infinite expansions. Revised and enlarged edition. D. Reidel Publishing Co., Dordrecht, 1974. 
\bibitem{3}
 Kim, D. S.; Kim, T. \emph{$r$-extended Lah-Bell numbers and polynomials associated with $r$-Lah numbers,} Proc. Jangjeon Math. Soc. \textbf{24} (2021), no. 1, 1-10.
\bibitem{4}
Kim, D. S.; Kim, T. \emph{Lah-Bell numbers and polynomials,} Proc. Jangjeon Math. Soc. \textbf{23} (2020), no. 4, 577-586.
\bibitem{5}
Kim, H. K. \emph{Degenerate Lah-Bell polynomials arising from degenerate Sheffer sequences,} Adv. Difference Equ. 2020, Paper No. 687, 16 pp. 
\bibitem{6}
Kim, D. S.; Kim, T. \emph{ A Note on a New Type of Degenerate Bernoulli Numbers,}  Russ. J. Math. Phys.\textbf{27} (2020), no. 2, 227-235.

\bibitem{7}
Kim, M.-S.; Kim, T. \emph{An explicit formula on the generalized Bernoulli number with order $n$,} Indian J. Pure Appl. Math. \textbf{31} (2000), no. 11, 1455-1461.
\bibitem{8}
Roman, S. \emph{The umbral calculus,} Pure and Applied Mathematics, 111. Academic Press, Inc. [Harcourt Brace Jovanovich, Publishers], New York, 1984. 
\end{thebibliography}
\end{document}